\documentclass[11pt]{article}
\usepackage{amsmath, amsthm, amssymb}    
\usepackage{bridges}			
\usepackage{graphicx}			
\usepackage[colorlinks=true, urlcolor=blue, citecolor=black, linkcolor=black]{hyperref}  
\usepackage{subcaption}			

\title{Quadratis Puzzles}

\author{Mario Guti\'errez\textsuperscript{1}, Hugo Parlier\textsuperscript{2} and Paul Turner\textsuperscript{3}
\vspace{10pt}\\
\textsuperscript{1}Mixed Reality Group, Logitech Europe; mgutierrez1@logitech.com\\
\textsuperscript{2}Department of Mathematics,
University of Luxembourg; hugo.parlier@uni.lu \\
\textsuperscript{3}Section de math\'ematiques, University of Geneva; paul.turner@unige.ch} 

\date{}					

\begin{document}

\maketitle

\thispagestyle{empty}

\begin{abstract}
We introduce a family of reconfiguration puzzles arising from ideas in geometry and topology. We present their construction from square-tiled shapes, discuss some of the underlying mathematics and describe how they are naturally associated to puzzle spaces which can be explored through visualization. 
\end{abstract}


\section*{Introduction}

There is a historic relationship between reconfiguration puzzles, such as the 15 puzzle or the Rubik's cube, and mathematics. In particular, there are many so-called ``sliding puzzles'' where pieces slide around a board (see for instance \cite{Hordern} for a collection of examples). In another direction, Weeks \cite{Weeks} took standard puzzles and transferred them onto a torus. 

In this paper we discuss a new family, Quadratis puzzles, which are inspired by topological and geometric constructions of surfaces. These puzzles also come with an associated ``space'' of possible configurations, that can be studied and even visualized. In contrast, the so-called Rubik's graph was studied for decades, resulting in a computation of its diameter \cite{Rokicki,RubikN}, but which is very far from being visualizable. The reader can also directly experience the puzzles for themselves by playing them \cite{website}.
\section*{The mathematics behind the puzzles}

\subsection*{Defining puzzles}
\begin{figure}[h]
\centering
\begin{minipage}[b]{0.25\textwidth} 
	\includegraphics[width=\textwidth]{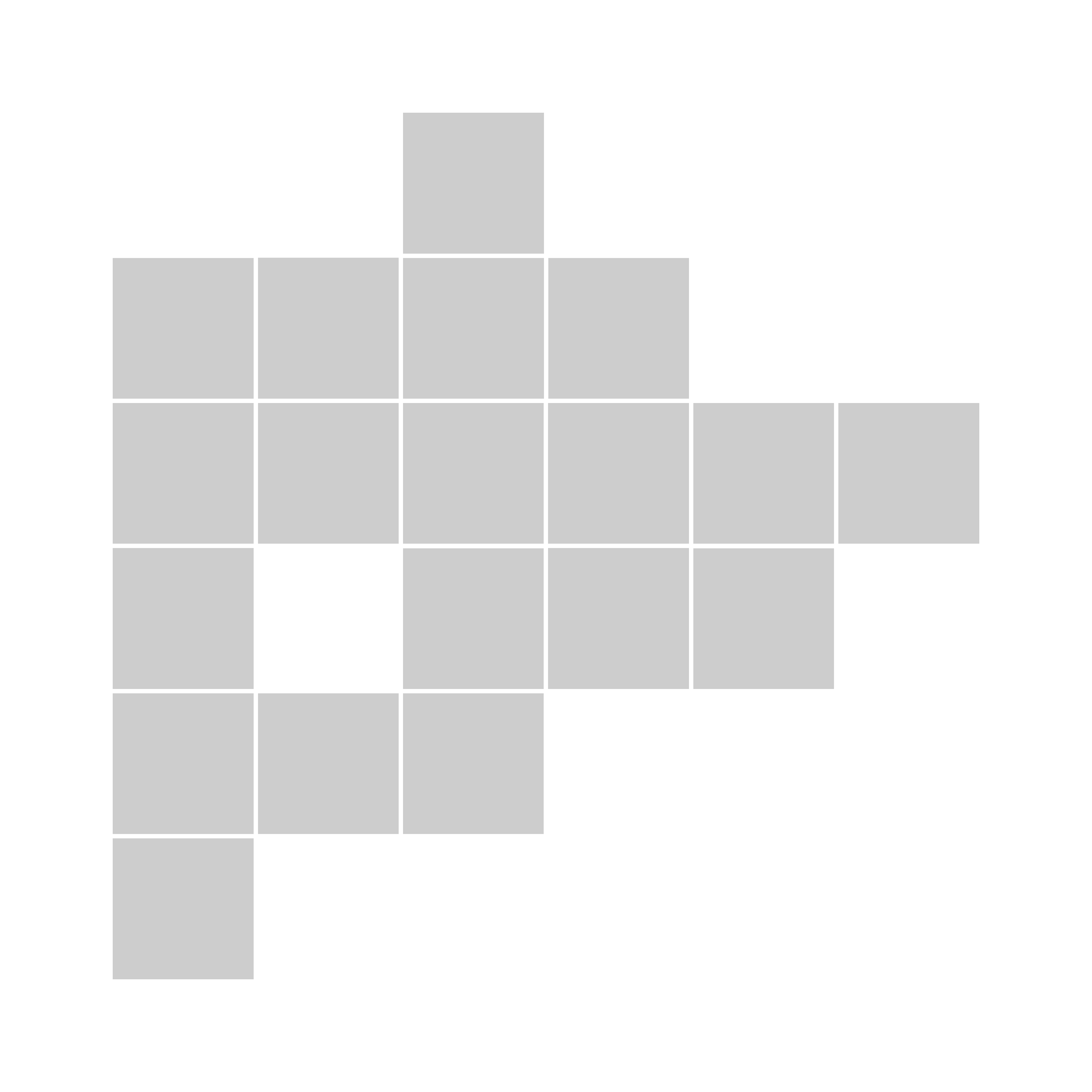}
        	\subcaption{A square-tiled shape} 
        \label{fig:boarda}
\end{minipage}	
\qquad\qquad
\begin{minipage}[b]{0.25\textwidth} 
	\includegraphics[width=\textwidth]{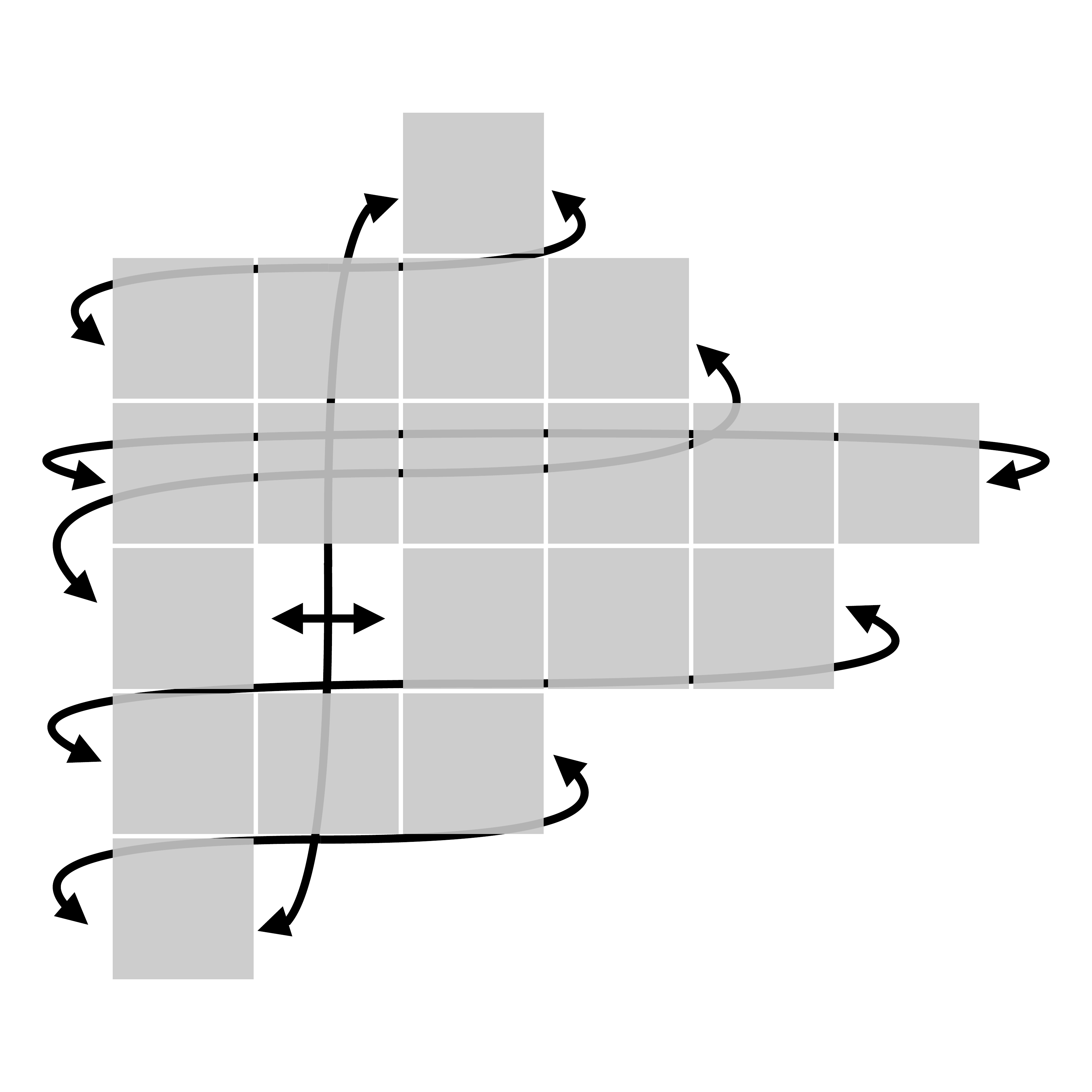}
        	\subcaption{Left-right pasting} 
       \label{fig:boardb}
\end{minipage}
\qquad\qquad
\begin{minipage}[b]{0.25\textwidth} 
	\includegraphics[width=\textwidth]{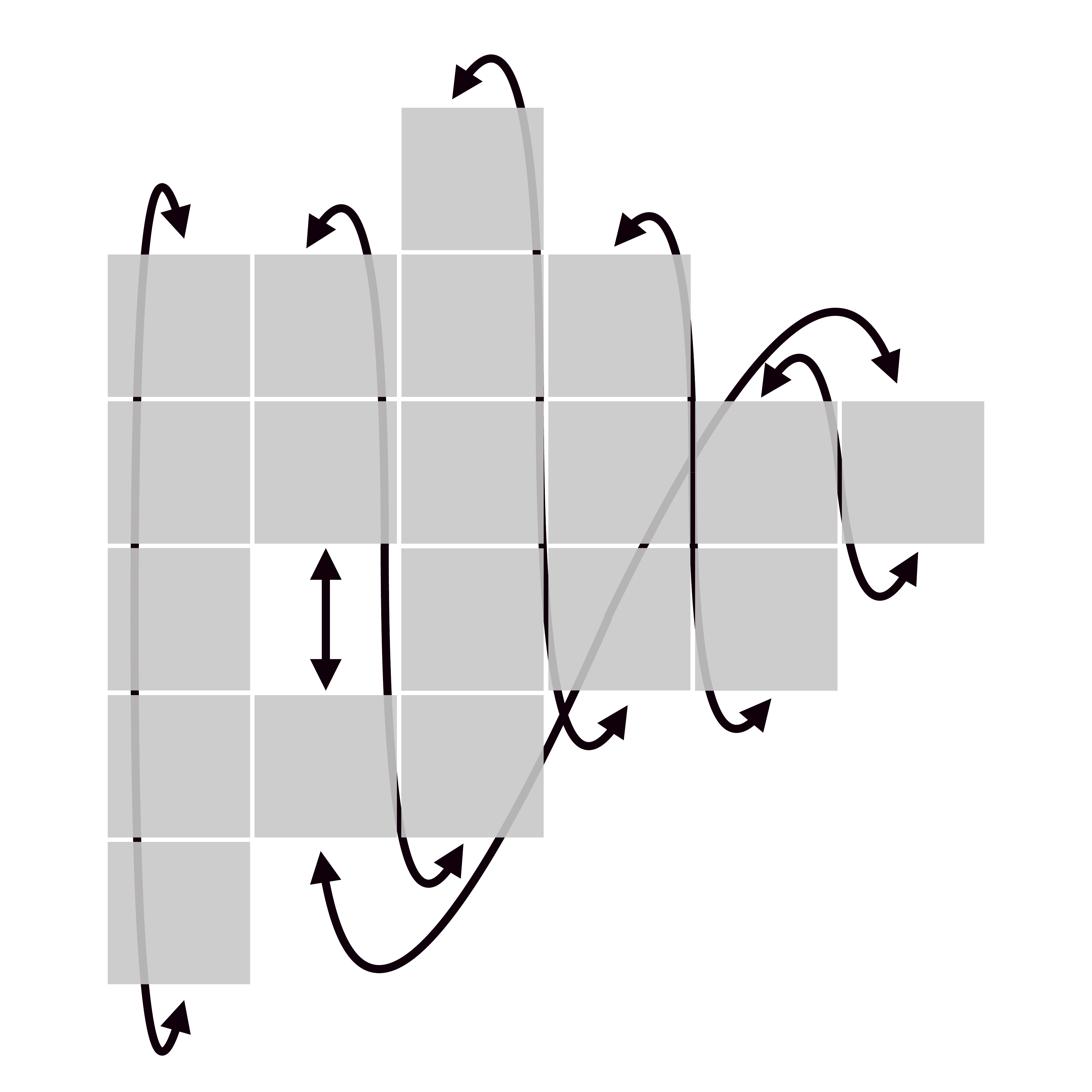}
        	\subcaption{Upper-lower pasting} 
        	\label{fig:boardc}
\end{minipage}
\caption{ }
 	\label{fig:board}
\end{figure}

To define a puzzle we must first describe what is meant by a {\em board} and a {\em color scheme}. A board is constructed by pasting together (unit) squares: we assemble collection of $N$ squares, each with upper, lower, left and right edges, by pasting  the edges together in such a way so that each left edge is pasted to a right edge, and each upper edge is pasted to a lower edge. This results in an orientable surface tiled by squares (see discussion below).
In practice, we often begin with a square-tiled shape that lives in the Euclidean plane, such as the one shown in Figure \ref{fig:boarda}. 

The non-pasted edges are then paired following the left-right or upper-lower conventions. One possible pasting of Figure \ref{fig:boarda}  is shown in Figure \ref{fig:boardb} (for left-right paring) and \ref{fig:boardc} (for upper-lower). 

\begin{figure}[h!tbp]
\centering
\begin{minipage}[b]{0.25\textwidth} 
	\includegraphics[width=\textwidth]{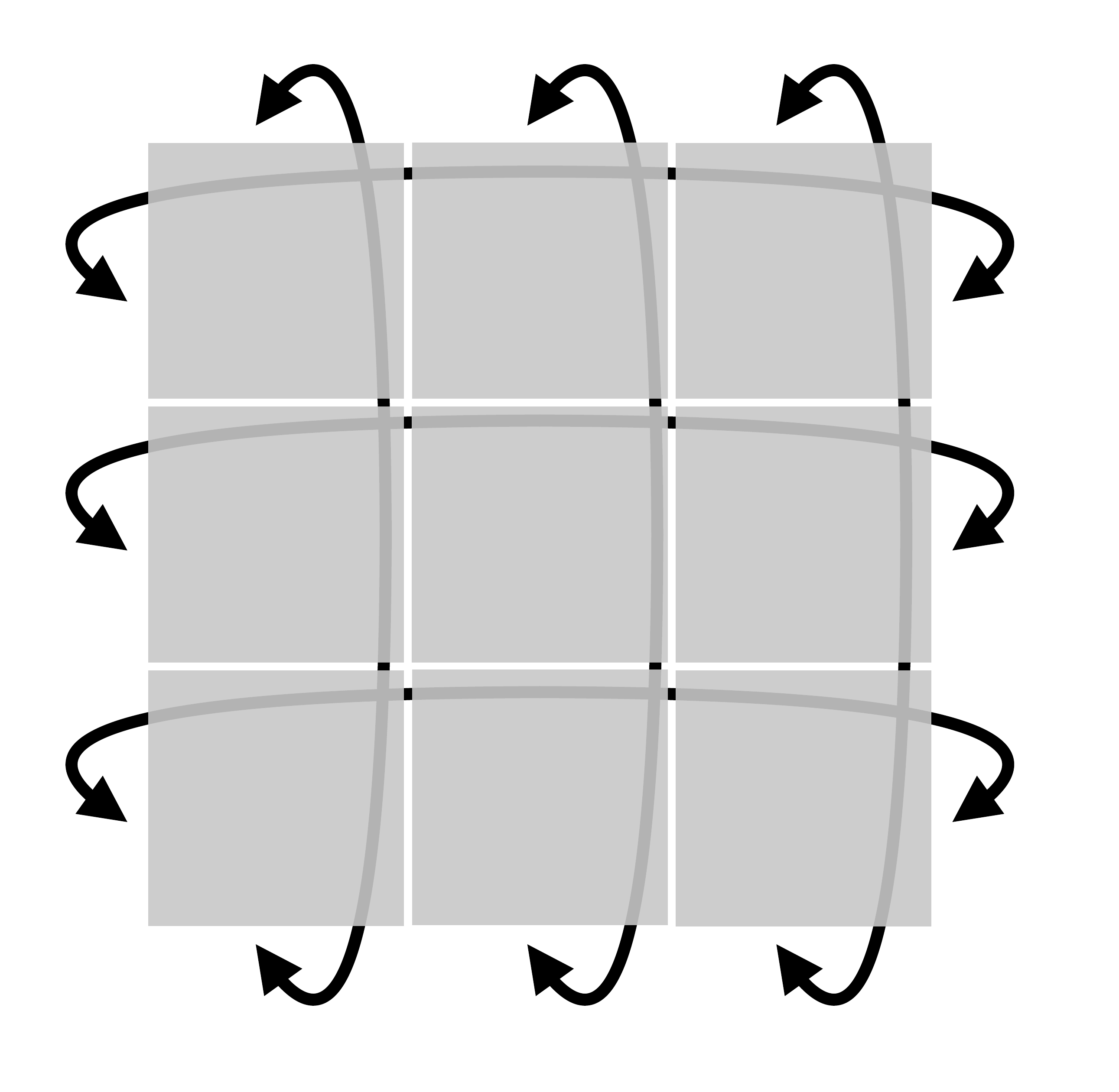}
        	\subcaption{} 
        \label{fig:squareboard}
\end{minipage}	
\qquad\qquad
\begin{minipage}[b]{0.25\textwidth} 
	\includegraphics[width=\textwidth]{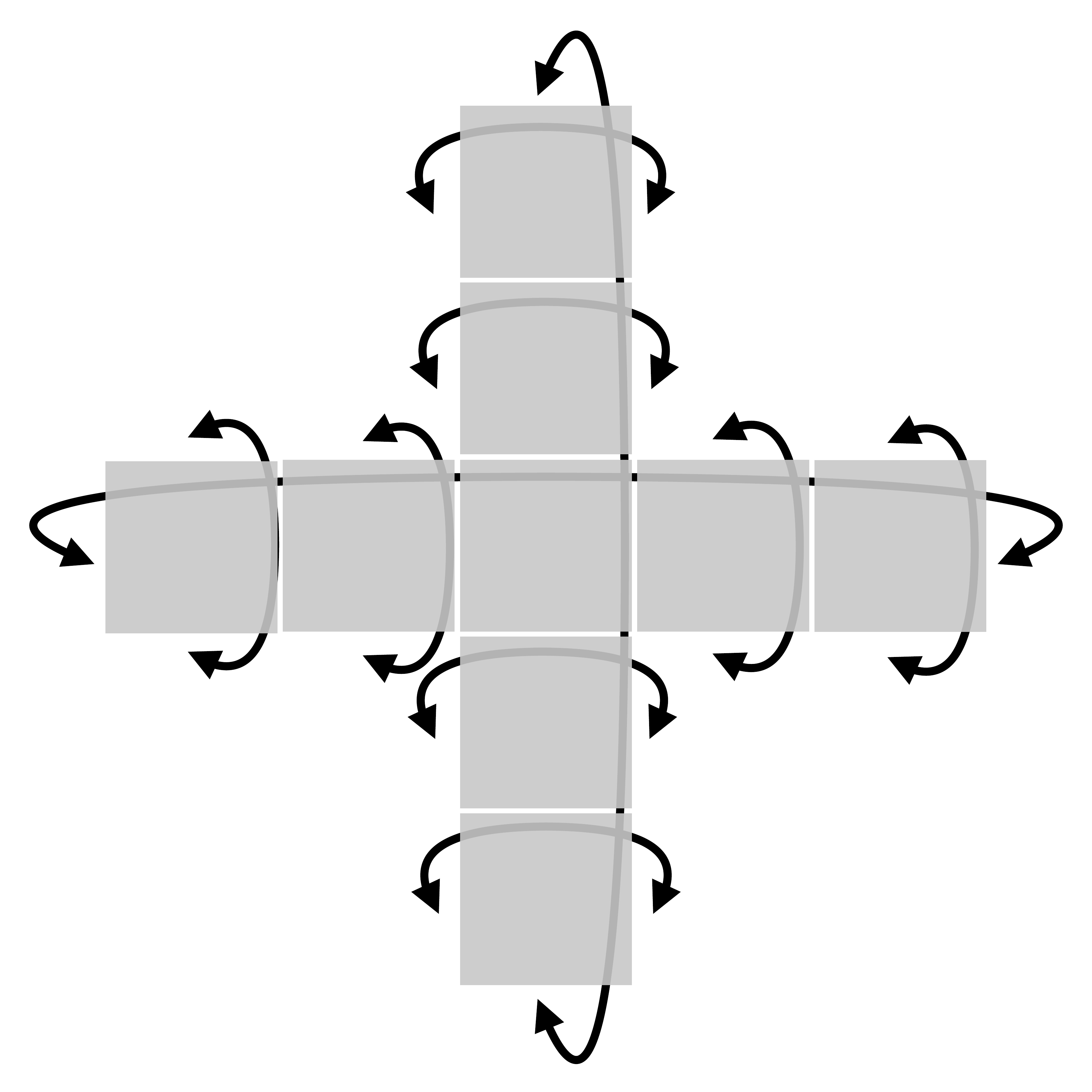}
        	\subcaption{} 
       \label{fig:crossboard}
\end{minipage}
\caption{Two square-tiled shapes and their standard pastings}
\label{fig:pastings}
\end{figure}

Given a square-tiled shape, there is an ``obvious'' way to paste the remaining sides together. Taking each horizontal line in turn, paste the first (ordering from left to right) free left side to the first available free right side on the same line. A similar procedure can be used to paste upper sides to lower sides in the same column. We'll refer to this pasting scheme as the {\em standard pasting}. In Figure \ref{fig:squareboard} and Figure \ref{fig:crossboard} we see the standard pastings of two different square-tiled shapes. Figure \ref{fig:boardb} and \ref{fig:boardc} present pasting that is not the standard pasting.

The squares of a board will be colored from a fixed palette with a predetermined number of squares of each given color. We refer to this information as a {\em color scheme}. Formally, though this is unnecessary for this article, we have a positive integer $K\leq N$ together with a $K$-tuple $(n_1,n_2,\ldots, n_K)$ such that $\sum_{n_i} = N$. Note that a color scheme does not determine how to color each square of a board, only the number of each color we must have.

\begin{figure}[h!tbp]
	\centering
	\includegraphics[width=0.8\textwidth]{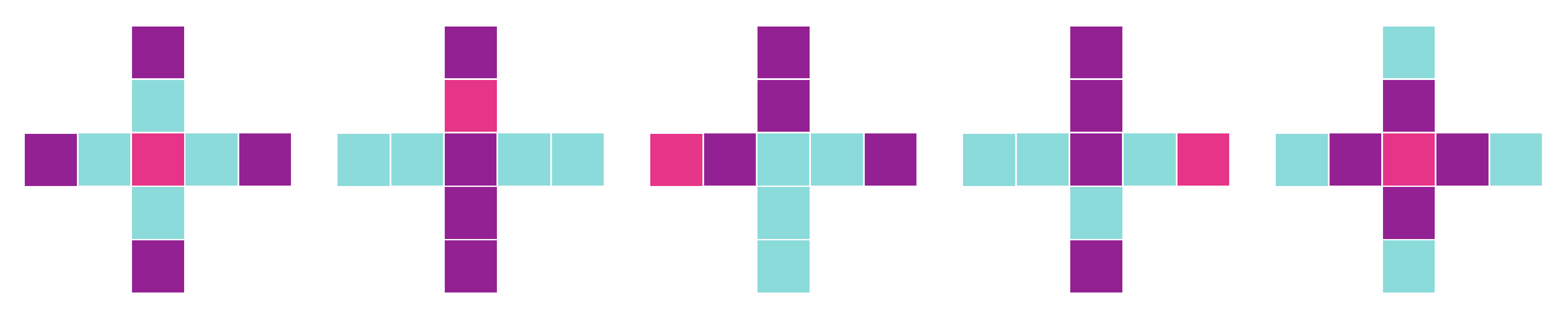}
	\caption{Five of the possible $630$ configurations with the given color scheme}
	\label{fig:configs}
\end{figure}

Given a color scheme a particular coloring of a board (of which there are many) will be referred to as a {\em configuration}. In Figure \ref{fig:configs} we see five different configurations of the board in Figure \ref{fig:crossboard} with the color scheme requiring 1 pink square, 4 blue squares and 4 purple squares.

\begin{figure}[h!tbp]
	\centering
	\includegraphics[width=0.7\textwidth]{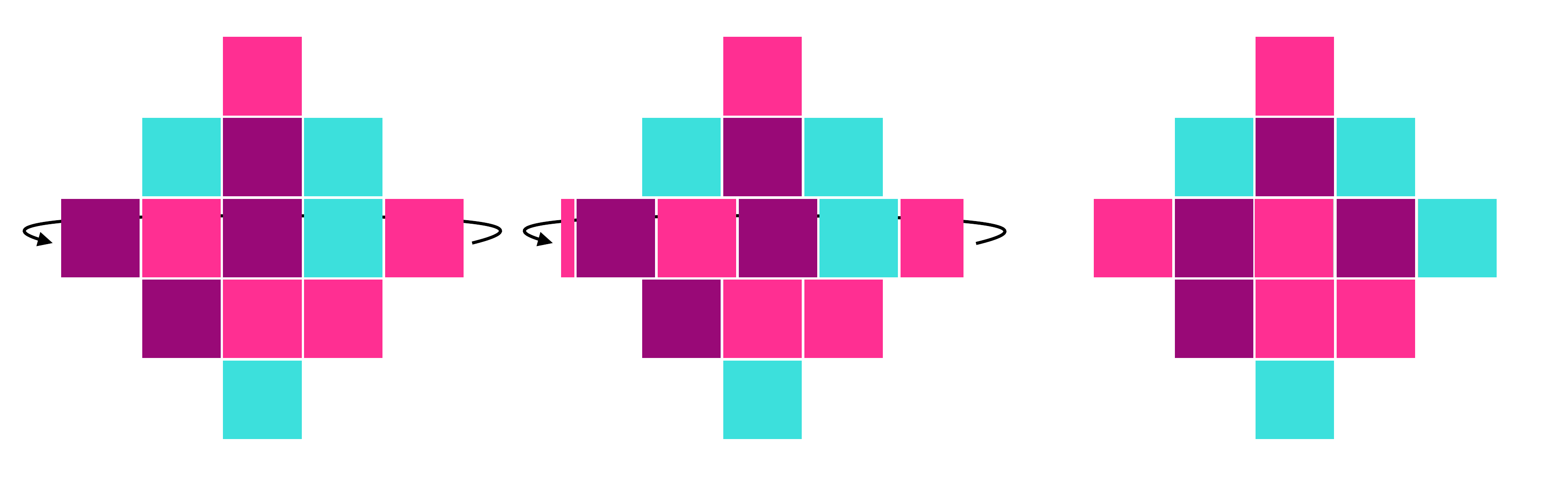}
	\caption{Making a move}
	\label{fig:move}
\end{figure}

To make a puzzle out of a board and a color scheme, we introduce a dynamic element, via {\em moves} which transform one configuration into another. A single square on a board always has $4$ (not necessarily distinct) neighboring squares: the ones above, below, to the right and to the left. ``Pushing'' a given square upwards, will have a knock-on effect on its upstairs neighbor, pushing it upwards as well. This in turn pushes its neighbor and so on and so forth until we end up back where we started. A coherent collection of squares are all displaced by one step upwards. A similar thing is seen if we push downwards, or to the left or to the right. We refer to these as {\em moves}. 
 Figure \ref{fig:move} illustrates a left to right horizontal move involving all squares of the middle row. This can be effected by ``pushing'' any square in that row towards the right. The reader is encouraged to visit \cite{website} to try this out for themselves.

Out of all of this we may now make a {\em puzzle}. Choose a particular configuration as the {\em home configuration}:  the aim will be to reproduce this configuration. Starting with the home configuration we can apply a number of random moves (in any direction) to obtain a new configuration. We refer to this as {\em shuffling}. The goal is to put the squares back into the home configuration, by using the moves described above. 
No distinction is made between two squares of the same color: if two squares of the same color are swapped in the process of solving, we do not care, the puzzle is still considered solved. Note that if a puzzle was shuffled using $m$ moves, it is entirely possible that it can be solved with fewer than $m$ moves. Pretty quickly, players tend to want to find the optimal number of moves to solve a puzzle, which brings us to the main mathematical object associated to a board equipped with a color scheme: its puzzle space. 

\subsection*{The configuration spaces of puzzles}

The {\em puzzle space} of a puzzle is a graph whose vertices are configurations and where two configurations are related by an edge if they can be related by a single move. We are only interested in configurations that can in principle be solved, so we simply ignore the other configurations. This amounts to only considering the component of the graph of all configurations which contains the puzzle home configuration. Thus, by definition, a puzzle space is a {\em connected} graph.

One could instead consider all possible colorings of a board with a prescribed set of colors, and then relate configurations as above. The graph obtained this way won't always be connected, and in fact determining the number and size of components seems to be an interesting problem. In \cite{Mathema}, there is an elementary proof of the fact that a standard 3 by 3 torus (see Figure \ref{fig:pastings} (a)) with all squares colored differently is disconnected.

\begin{figure}[h!tbp]
	\centering
	\includegraphics[width=.5\linewidth]{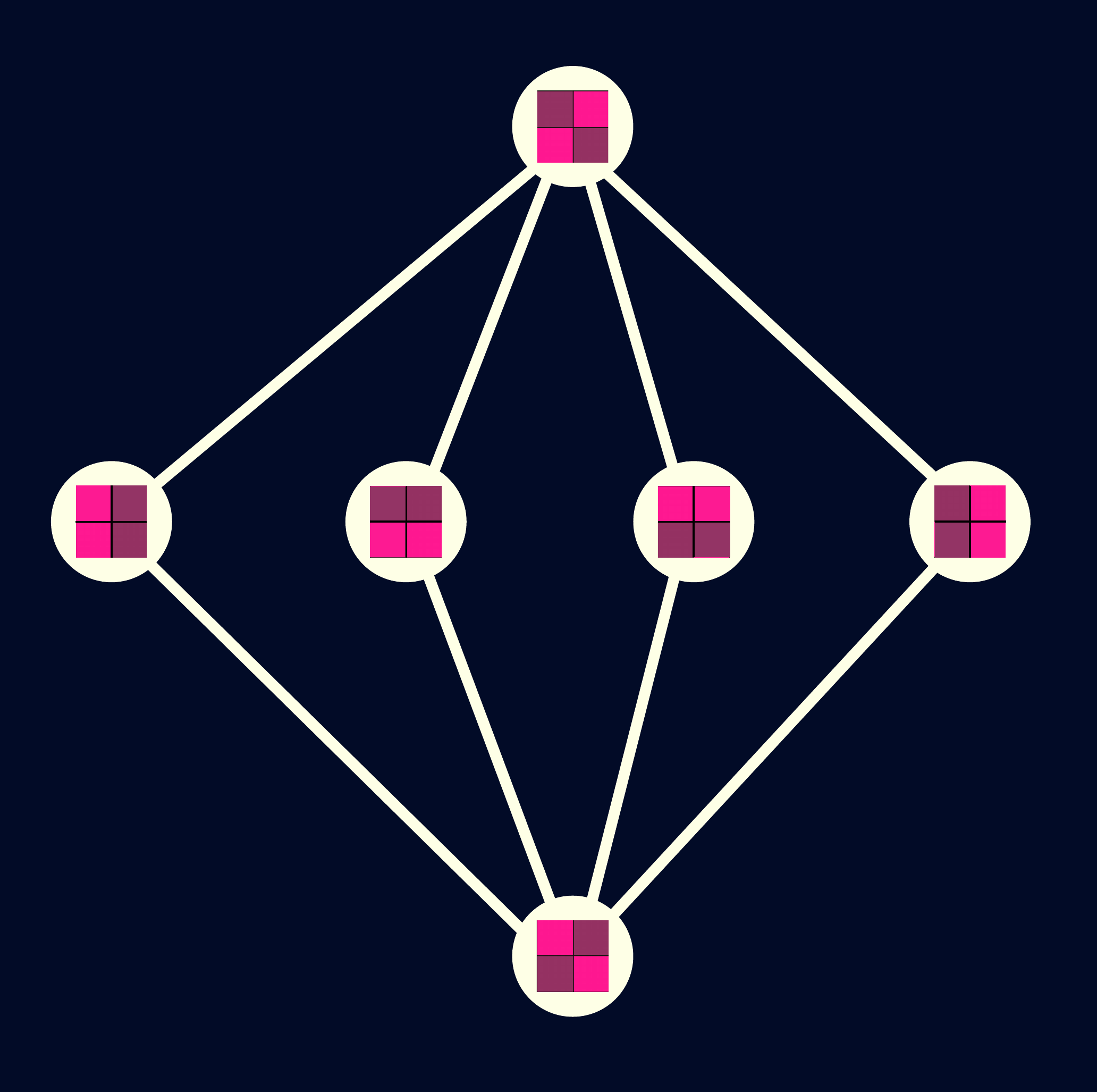}
	\caption{A simple puzzle space}
	\label{fig:simple}
\end{figure}

In Figure \ref{fig:simple} we see a simple example which can be worked out by hand. The board is a two-by-two square with the standard pasting, and we color the squares with two colors and two squares of each color. Any of the 6 configurations shown may be used as the home configuration.

\begin{figure}[h!tbp]
\centering
\begin{minipage}[b]{0.49\textwidth} 
	\includegraphics[width=\textwidth]{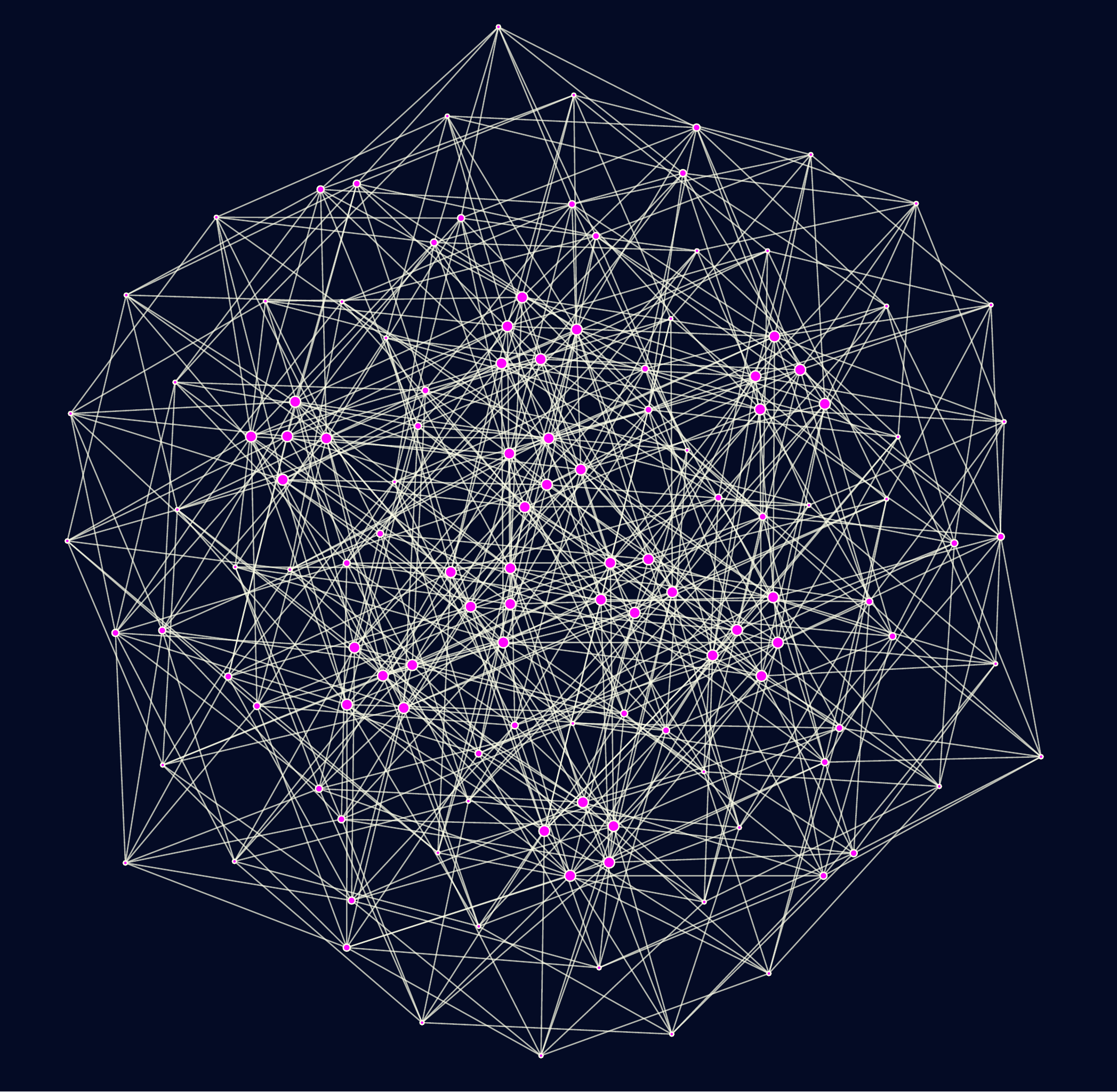}
		\subcaption{}
   \label{fig:three}
        
\end{minipage}
~ 
\begin{minipage}[b]{0.49\textwidth} 
	\includegraphics[width=\textwidth]{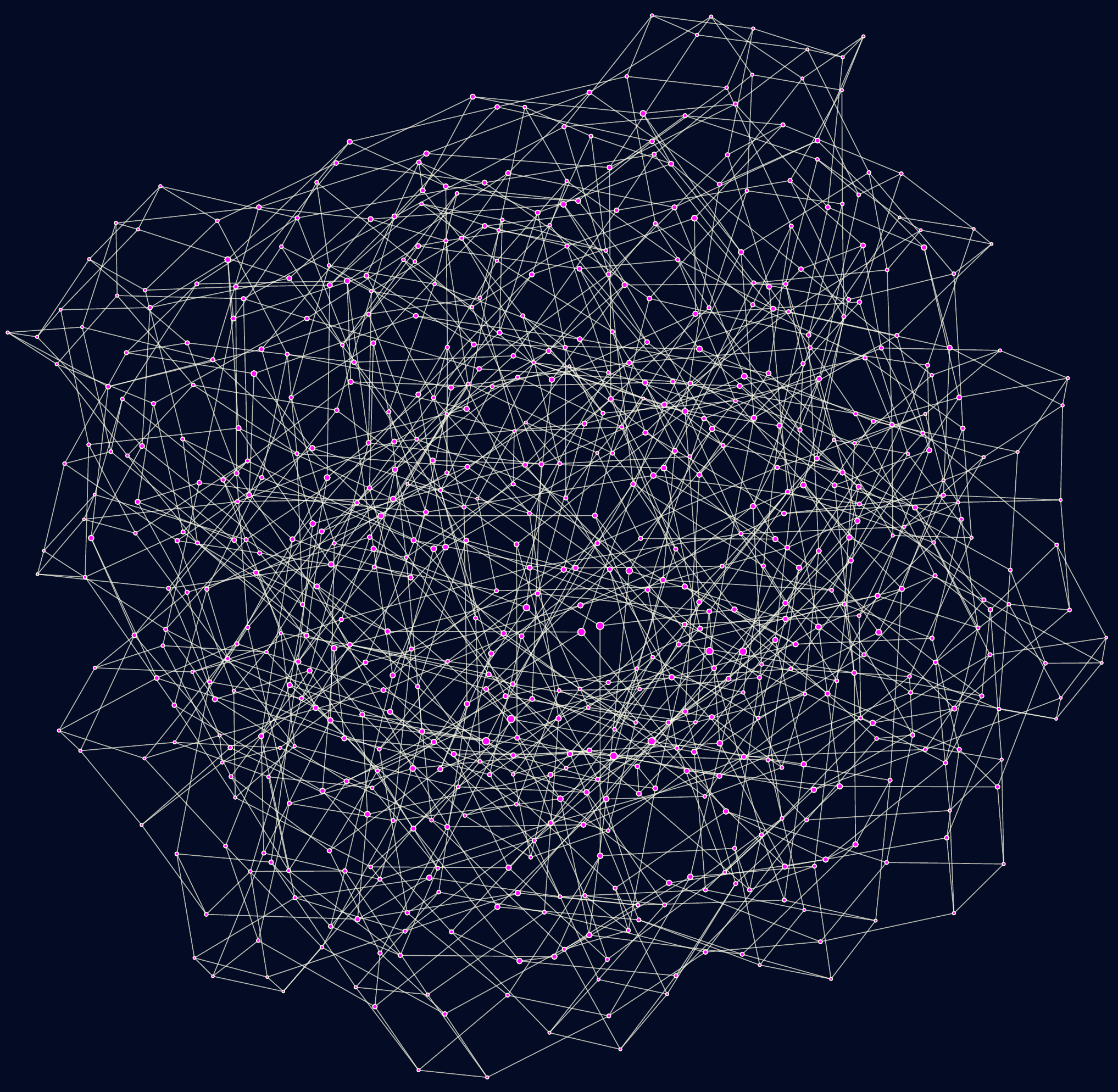}
        	\subcaption{}
		\label{fig:crossgraph}
       
\end{minipage}
\caption{Puzzles spaces display intriguing features}
\end{figure}

As the complexity of the puzzle increases, the graphs become more complicated and enjoy richer geometries. For instance, Figure \ref{fig:three} is the graph for a three-by-three square board with standard pasting and two colors, 3 of one color and 6 of the other. As all configurations are related by moves in this example, the choice of home configuration is again unimportant. The figure represents one of many possible planar representations of the graph, and was obtained using standard graph visualisation software (Gephi) and a force-directed graph layout algorithm.

\subsection*{Relationships to moduli theory}

By standard surface topology techniques the pasting procedure described above results in an orientable surface equipped with a decomposition into squares. In the examples shown in Figure \ref{fig:board} and Figure \ref{fig:pastings}, the underlying surface is connected, but that is not necessarily a requirement. Interestingly, one can show that if one randomly pastes the squares together, the result is often connected. In fact, the probability of obtaining a connected surface goes to $1$ as the number of squares goes to infinity \cite{Shrestha-Yang,Shrestha}.

The board given by Figure \ref{fig:crossboard}  is a genus $2$ surface (so topologically a double-torus). This can be shown by a standard Euler characteristic computation and the classification of surfaces. By using the Euler characteristic again, it can be proved that the board from Figure \ref{fig:board} is a genus $6$ surface. 

Beyond the topology and the combinatorics of the decomposition into squares, the boards are surfaces with a geometric structure. In fact, they are important objects in the theory of moduli spaces, and are generally called square-tiled translation surfaces. (The word ``translation'' comes from the fact that side identifications are made by Euclidean translations.) These surfaces have all sorts of interesting properties, including being dense in the space of all translation surfaces (up to scaling). Some of them are quite exotic and have fascinating dynamical and geometric properties \cite{Schmithusen} and there has been a lot of research about square-tiled and translation surfaces in general, a partial survey can be found in \cite{Forni-Matheus}. As an example, a very natural question is:  how many square-tiled translation surfaces can be made with $N$ squares? It turns out to be very difficult and is related to computing volumes of moduli space (see \cite{Eskin-Okounkov} for partial results). In general the answer is not known, making the question of enumerating the number of possible puzzle boards with a given number of squares, an open and difficult problem.

For these reasons, the puzzle spaces we introduce above can be seen as ``toy'' moduli spaces, and can be studied using a similar approach. How many vertices are there? How do you determine the distance between two configurations? More generally, what do these spaces ``look'' like? 

These questions, as for moduli spaces, can range from specific to more general, from elementary to more sophisticated. For instance, what are their diameters? Whereas a very limited number of small examples can be worked out, finding a general answer to this type of question very quickly turns into serious and seemingly difficult research questions.

\section*{Visualizing puzzle spaces}

There are two desirable properties we might ask of a puzzle, namely that it is {\em interesting}  and {\em visualizable}. These are rather loose concepts:  interesting means non-trivial yet doable in the sense that a player can, with some thought, develop strategies which help to provide a solution; visualizable means that the puzzle space is of a size that may be displayed so that the vertices and edges remain discernible. These two constraints pull in different directions with increasing interest often resulting in graphs that are too large. 

\subsection*{The first examples: chroma-squares}
A particularly simple class of puzzles was introduced in \cite{Mathema} under the name {\em chroma-squares}:  the board is a square $N\times N$ array and the color scheme has $N$ colors with  $N$ squares of each. The pasting is the standard one, meaning the underlying surface is a torus.

\begin{figure}[h!tbp]
	\centering
	\includegraphics[width=\linewidth]{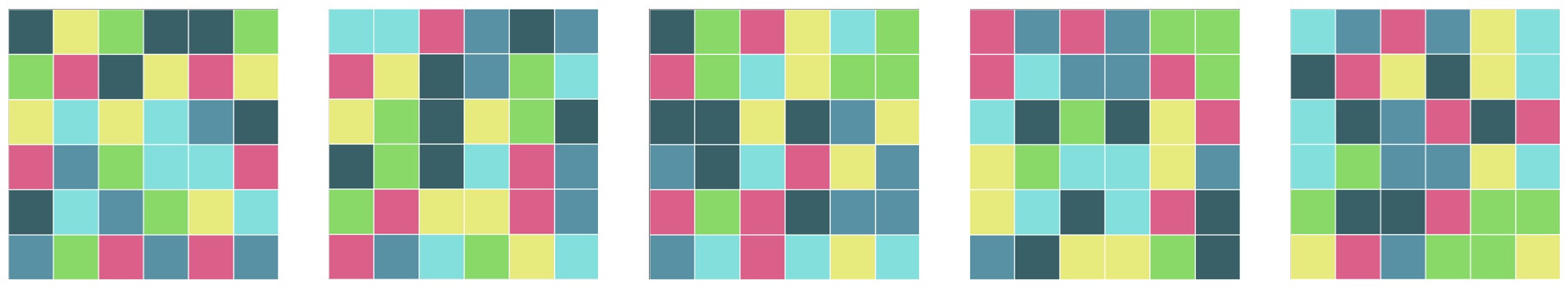}
	\caption{Five of the approximately $2.7 \times 10^{24}$ configurations of size $6 \times 6$, images from \cite{Mathema}}
	\label{fig:1}
\end{figure}

There are some rather obvious choices to make for home configurations, namely lines or columns of uniform color. With essentially no mathematical background it is possible to explain that any two configurations can be linked by a sequence of moves (see \cite{Mathema}), which shows that the graph of {\em all} configurations is connected.

With 1680 vertices and around 10,000 edges, the case $N=3$ yields a puzzle that is both {\em interesting} and {\em visualizable}. For $N=2$ the puzzle is too trivial and the graph too simple (it is given in Figure \ref{fig:simple} above) and for $N>3$ the graph is too large (although the difficultly of solving the puzzle does not really increase). 

\subsection*{Beyond tori}
As we saw above, the underlying surface of a puzzle is an orientable surface. Higher genus examples also give a rich source of interesting and visualizable puzzle spaces. Recall that the board shown in Figure \ref{fig:crossboard} has genus two. With color scheme indicated in Figure \ref{fig:configs} it is a relatively easy puzzle (for any chosen home configuration), but not completely trivial, and it possesses a visualizable puzzle graph as shown in Figure \ref{fig:crossgraph}.




\subsection*{Exploring large puzzles}\label{sec:large}

\begin{figure}[h!tbp]
\centering
\begin{minipage}[b]{0.3\textwidth} 
	\includegraphics[width=\textwidth]{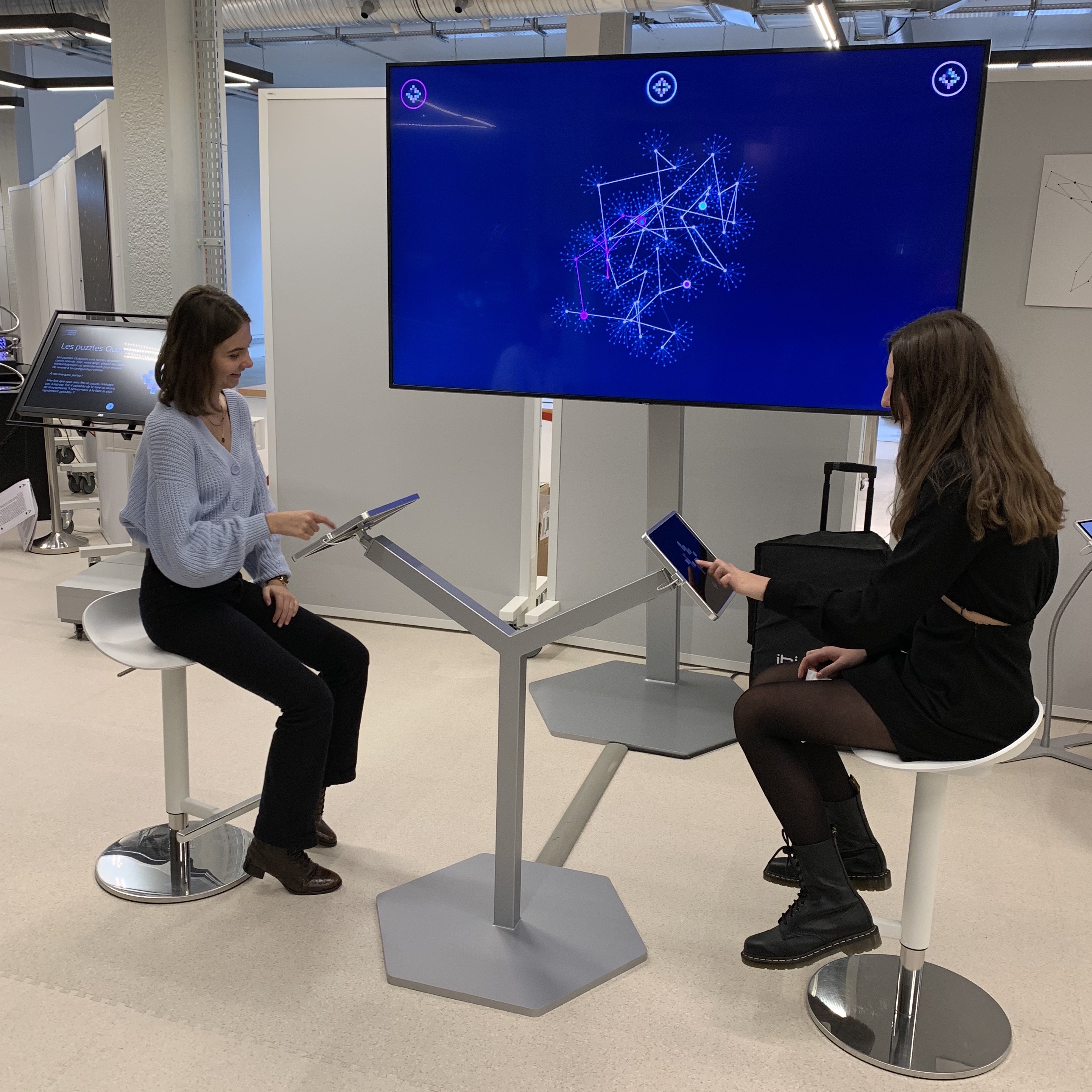}
		\subcaption{{\em Exploratis} installation at the Luxembourg Science Center}
    	\label{fig:lux}
        
\end{minipage}
~ 
\qquad\qquad
\begin{minipage}[b]{0.3\textwidth} 
	\includegraphics[width=\textwidth]{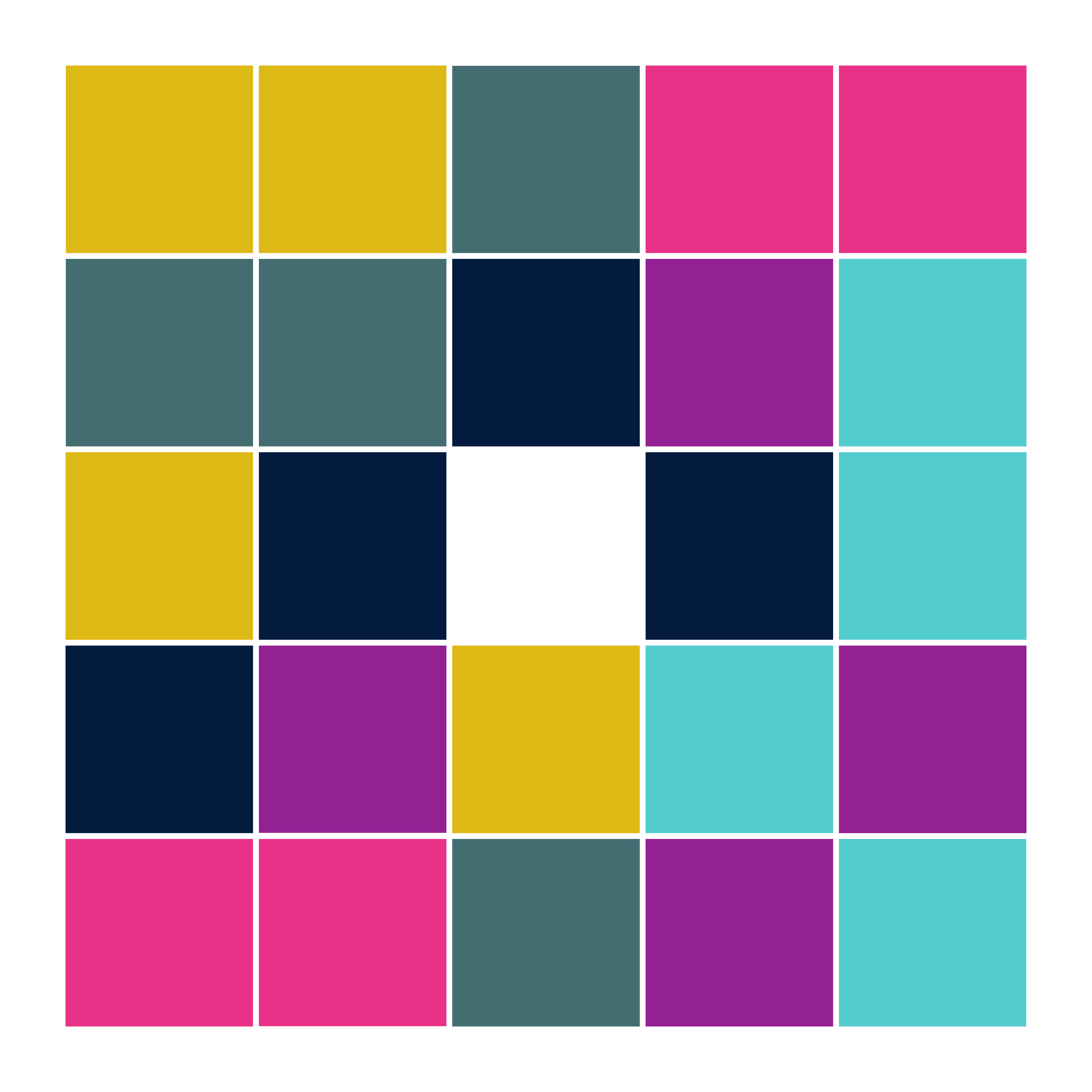}
        	\subcaption{A puzzle which is too big to be explored entirely}
		\label{fig:hole}
       
\end{minipage}
\caption{}
\end{figure}

While many  puzzles are both interesting and visualizable, sometimes we may want to partially visualize a configuration space of a puzzle complicated enough so the entire space is beyond visualization (either for computational or presentational reasons). The puzzle shown in Figure \ref{fig:hole} (in which the board is a $5\times 5$ grid missing its middle square) is interesting because the usual strategies do not resolve a potential problem of switching at the end. There is something new here, which means it is worth retaining as an example despite its large configuration space (which has has over 3000 billion billion points, roughly the number of square millimeters needed to cover the entirety of the Grand Duchy of Luxembourg). 

One solution is to dynamically develop the graph as the player walks through different configurations: if a new configuration is reached, add it to the graph along with the edge used to reach it. This is the basis of the Science Museum installation {\em Exploratis} developed in conjunction with the Luxembourg Science Center in which an iPad driven puzzle is linked to a large screen displaying an ever larger portion of the configuration space (see Figure \ref{fig:lux}). Here again, a type of repulsion based algorithm is used to visualize the graph. Players can play individually, race or collaborate. This time, because the graph changes with every move, the visuals have a wow-effect. 
\begin{figure}[h!tbp]
\centering
\begin{minipage}[b]{0.28\textwidth} 
	\includegraphics[width=\textwidth]{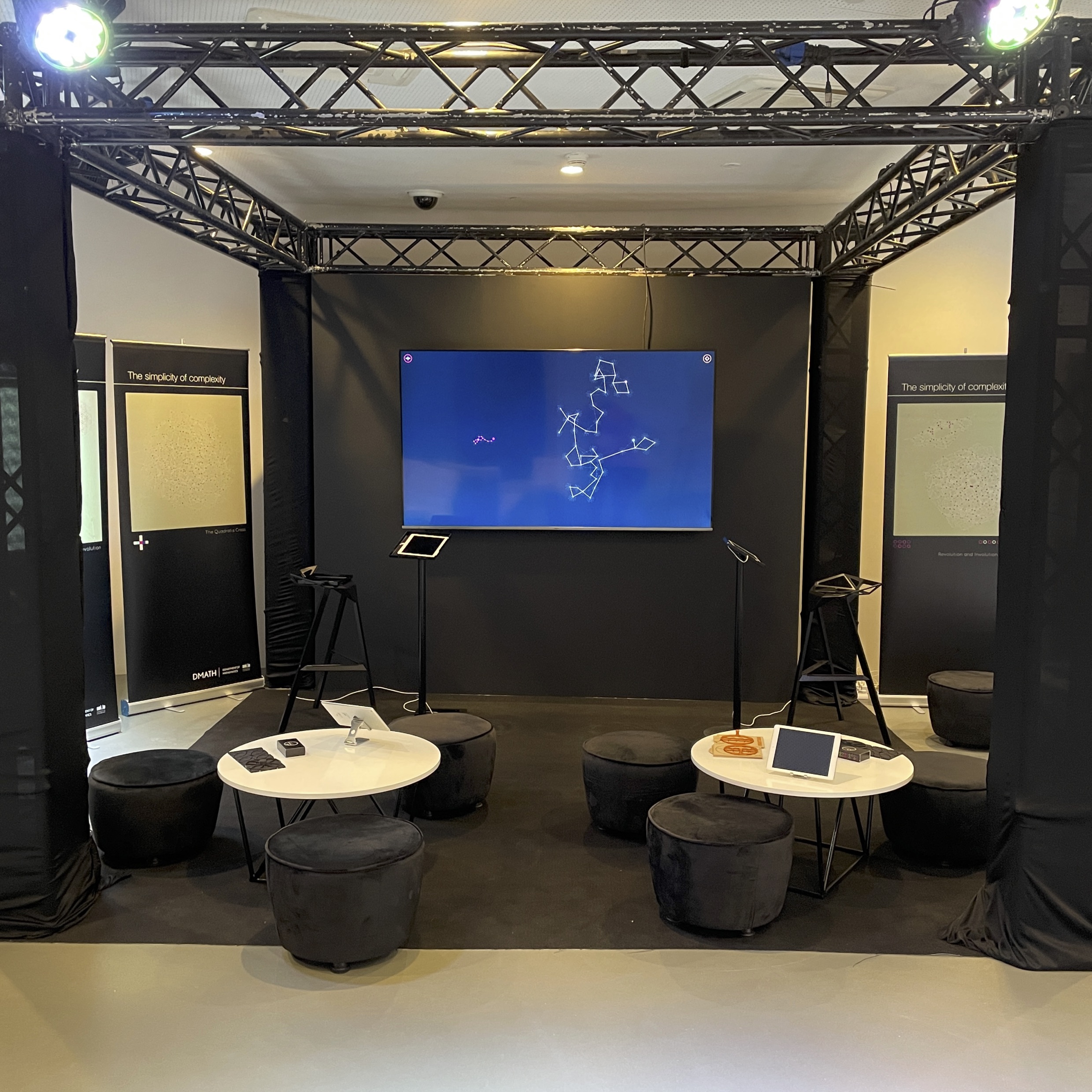}
        	\subcaption{} 
        	\label{fig:2b}
\end{minipage}
\qquad
\begin{minipage}[b]{0.28\textwidth} 
	\includegraphics[width=\textwidth]{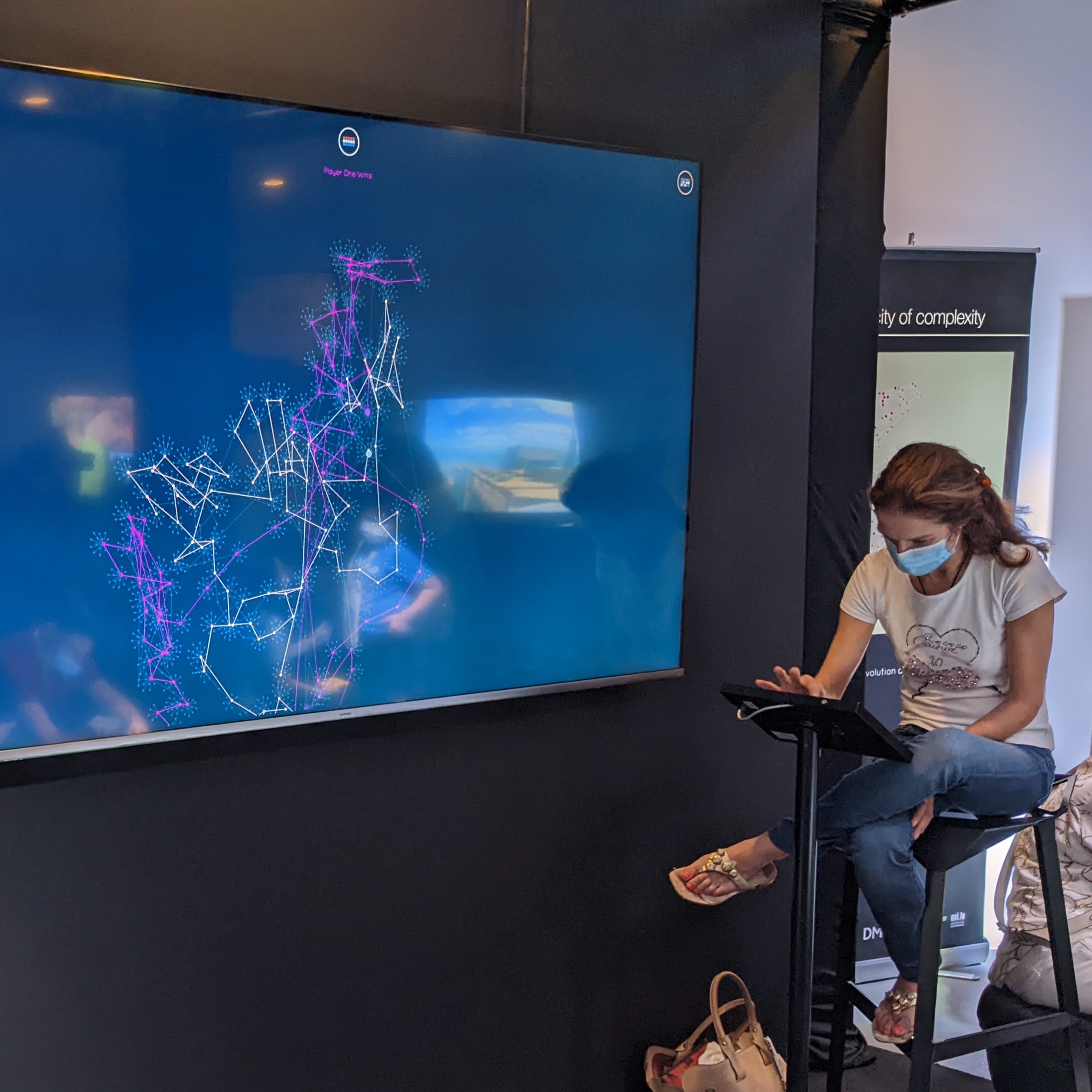}
        	\subcaption{} 
        	\label{fig:2a}
\end{minipage}
~ 
\qquad
\begin{minipage}[b]{0.28\textwidth} 
	\includegraphics[width=\textwidth]{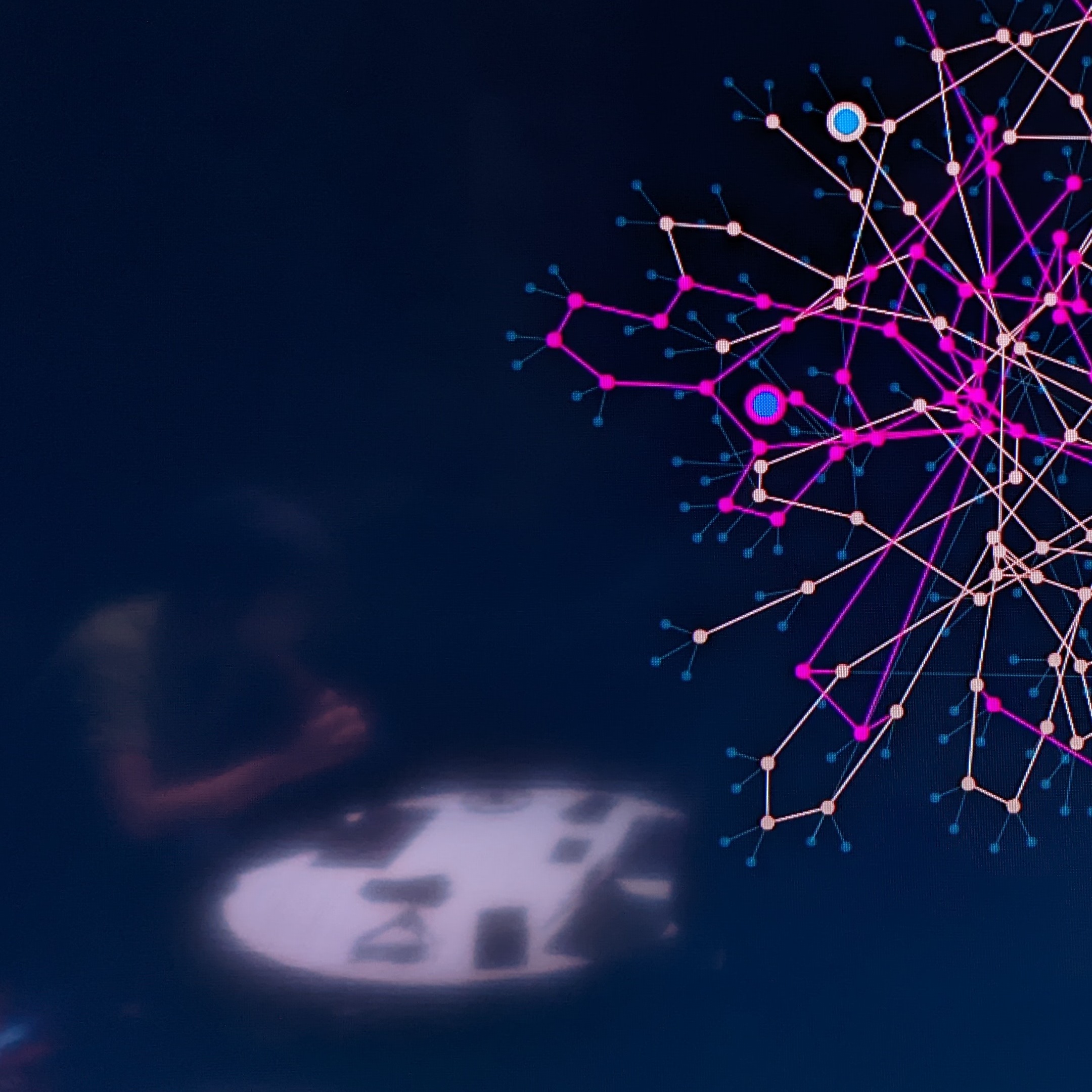}
        	\subcaption{} 
        	\label{fig:2b}
\end{minipage}
\caption{ {\em Exploratis} installation at Expo Dubai (Luxembourg Pavilion) }
\label{fig:dubai}
\end{figure}
A similar installation was made available to visitors of the Luxembourg Pavilion at the World Expo in Dubai. Hundreds of visitors each day were able to test the puzzles and explore the graphs. In each case, the puzzles are tailored for the event. While many of the puzzles that we use are the same as the one in the available app, we can create new puzzles, such as one representing the flag of Luxembourg for the pavilion, for special installations.

\section*{Conclusion}
Although only recently made public, we have already had several opportunities to test the Quadratis puzzles with a variety of different player profiles, and we get continual feedback from the Science Center. Something that was not obvious to us at first, and that required a lot of testing, is that the puzzles can be scaled in difficulty, ranging from completely obvious to (apparently) beyond the reach of human capability (and with basically everything in between). By careful planning, a sequence of puzzles can be presented that allows players of any level to learn as they progress. The lessons learnt, along with tricks and techniques uncovered through play, can be subsequently applied to increasingly challenging situations. 


So why is this interesting? The puzzles can be appreciated without any scientific background, and yet the surrounding mathematics gives a player a glimpse at what it is like to explore unchartered mathematical territory. In light of the feedback we've received so far, we see this is a starting point for a multi-faceted adventure. What is particularly exciting  is to see the diversity of players (background, origin, age...) that have experienced the joy of simply playing these puzzles. As we create more and more puzzles, we are ourselves discovering new and interesting phenomena. The exploration of the mathematics, particularly the study of the associated graphs, has only just started. In a different direction, it would be fascinating to know more about what and how players learn while progressing though a sequence of puzzles. 
 
Quadratis puzzles are regarded as simple fun by some, but the fact that they can be used to illustrate mathematics and the nature of mathematical research remains, for us, a central motivation.
 
\section*{Acknowledgements}
This article is part of a more general long standing project aiming to illustrate mathematical research in different ways, and we've received a lot of support along the way.

The original Mathema project was co-financed by the Swiss National Science Foundation and by TomBooks. In particular, we thank Thomas Steinmann for his personal support.

The exhibitions where we have presented Quadratis and Exploratis are brought alive by mediators so many thanks to all of you who have participated in science festivals. There are too many of you to name, but a particular thank you goes to Bruno Teheux. Not only was he co-organiser of many of the events but his feedback and encouragement has been invaluable. More generally, the University of Luxembourg, and in particular the Department of Mathematics of the Faculty of Science, Technology and Medicine, has provided logistical and financial aid. Thank you in particular to Giovanni Peccati. We also thank Sorbonne Abu Dhabi and the University of Geneva for help with certain exhibitions and to the Luxembourg Pavillion at the Expo in Dubai for hosting us.

The Exploratis station had the support of many partners. First of all, thank you to the Luxembourg Science Center for hosting us, and to Julien Meyer for his help in making this vision a reality. Exploratis was supported by a PSP Classic grant from the Luxembourg National Science Foundation, who has also supported the project indirectly in many ways (Researchers' Days and Science Festivals). Exploratis was also co-developed with Reyna Ju\'arez who went far and beyond the call of duty. 

    
{\setlength{\baselineskip}{13pt} 
\raggedright				

} 
   
\end{document}